\input amstex
\documentstyle{amsppt}
\magnification=\magstep1

\hoffset1 true pc
\voffset2 true pc
\hsize36 true pc
\vsize50 true pc

\tolerance=2000 
\define\m1{^{-1}}
\define\ov1{\overline}
\def\gp#1{\langle#1\rangle}
\def\ul2#1{\underline{\underline{#1}}}

\def\rad{ \operatorname{rad \,}}

\TagsOnRight

\topmatter
\title 
On a filtered multiplicative basis of group algebras
\endtitle

\author
Victor Bovdi
\endauthor

\leftheadtext\nofrills{Victor Bovdi}
\rightheadtext\nofrills{
On a filtered multiplicative basis of group algebras
}

\abstract
Let $K$ be a field  of characteristic $p$ and $G$ a nonabelian 
metacyclic  finite $p$-group. 
We give an explicit list of all metacyclic $p$-groups $G$, such
that the group algebra $KG$ over a  field of characteristic $p$ 
has a  filtered multiplicative $K$-basis. We also present an 
example of a non-metacyclic $2$-group $G$, such that the group
algebra $KG$ over any field of characteristic $2$ has a filtered 
multiplicative $K$-basis.

\endabstract

\thanks 
The first author's research was supported  by the Hungarian 
National Foundation for Scientific Research No.T 025029  
and by FAPESP Brazil (proc. 97/05920-6).
\endthanks 

\address
\hskip-\parindent
Institute of Mathematics and Informatics
\newline
Lajos Kossuth University
\newline
H-410 Debrecen, P.O. Box 12
\newline
Hungary
\newline
vbovdi\@math.klte.hu
\endaddress

\endtopmatter
\document

\subhead
1.Introduction
\endsubhead
Let $A$ be a finite-dimensional algebra over a field $K$ and let $B$ be a
$K$-basis of $A$.  Suppose that $B$ has the following properties:
\itemitem{1.} if $b_1,b_2\in B$ then either $b_1b_2=0$ or $b_1b_2\in B$;
\itemitem{2.}  $B\cap \rad(A)$ is a $K$-basis for $\rad(A)$, where 
$\rad(A)$ denotes  the Jacobson radical of $A$.

Then $B$ is called a {\it  filtered multiplicative $K$-basis} of $A$. 

The filtered multiplicative $K$-basis   arises  in the theory of
representations of algebras and was introduced first by H.~Kupisch \cite{5}. 
In \cite{1} R.~Bautista,
P.~Gabriel, A.~Roiter and L.~Salmeron  proved that if there are only 
finitely many isomorphism classes of indecomposable $A$-modules over an algebraically 
closed field $K$,  then $A$  has   a filtered multiplicative $K$-basis. Note 
that by Higman's theorem  the group algebra $KG$ over a   field of 
characteristic $p$ has only finitely many isomorphism classes of
indecomposable $KG$-modules if  and only if all the Sylow $p$-subgroups of $G$
are cyclic.

Here  we study  the following question from \cite{1}: {\it When does 
exist   a filtered multiplicative $K$-basis in the group algebra $KG$}?

Let $G$ be a finite abelian $p$-group. Then  
$G=\gp{a_1}\times \gp{a_2}\times \cdots \times \gp{a_s}$ is the direct product
of  cyclic  groups $\gp{a_i}$ of order  $q_i$,   the set
$$
B=\{(a_1-1)^{n_1}(a_2-1)^{n_2}\cdots (a_s-1)^{n_s}\,\,\,\, \mid\,\,\,\,   0\leq n_i< q_i\}
$$
is  a filtered multiplicative $K$-basis of the group algebra $KG$ over the
field $K$ of characteristic $p$.

Moreover, if $KG_1$ and $KG_2$ have  filtered multiplicative $K$-bases,
which we call $B_1$ and $B_2$ respectively, 
then $B_1\times B_2$  is a filtered multiplicative $K$-basis  of the group algebra
$K[G_1\times G_2]$. 

P.~Landrock and G.O.~Michler \cite{6} proved that the group algebra of the 
smallest Janko group over a field of  characteristic $2$ does not have a  filtered
multiplicative $K$-basis.  

L.~Paris \cite{7} gave  examples of  group algebras  $KG$, which have  no filtered
multiplicative $K$-bases. He also showed  that if $K$ is a field of characteristic
$2$ and either a) $G$ is a quaternion group of order $8$ and also $K$  contains
a  primitive cube root of the unity or b) $G$ is a dihedral $2$-group,  then
$KG$  has a filtered  multiplicative $K$-basis. We shall show that  for the
class of all metacyclic $p$-groups the groups mentioned in the items a) and b) are
exactly those  for which a  multiplicative $K$-basis exists.

We also present an example of a non-metacyclic $2$-group $G$, such that the 
group algebra $KG$ over any field of characteristic $2$ has a filtered 
multiplicative $K$-basis.

\subhead
2. Preliminary remarks and notations
\endsubhead
Let $B$ be  a filtered multiplicative $K$-basis  in 
a finite-dimensional $K$-algebra $A$. In the proof of the  main
result we use the following simple properties  of $B$:
\itemitem{}(I) $B\cap rad(A)^n$ is a $K$-basis of 
$rad(A)^n$ for all $n\geq 1$.

Indeed, by the definition of a  basis, $B\cap rad(A)$ is a $K$-basis of $rad(A)$ and the
subset  $B\cap rad(A)^n$ is linearly independent over $K$. Since the set of
products $b_1b_2\cdots b_n$ with $b_i\in B$ is a generator system for
$rad(A)^n$ and each such product is  either $0$ or belongs to $B\cap rad(A)^n$,
we conclude that $B\cap rad(A)^n$ is a $K$-basis of $rad(A)^n$. 
\itemitem{}(II) if $u,v\in B\setminus rad(A)^k$ and 
$u\equiv v \pmod{rad(A)^k}$ then $u=v$. 

Indeed, if  $u-v=\sum_{w\in B\cap rad(A)^k}\lambda_ww$  with
$\lambda_w\in K$, then by the  linearly  independency of the basis elements we
conclude that $\lambda_w=0$ and therefore $u=v$.

Recall that  the Frattini subalgebra $\Phi(A)$ of $A$ is defined as the
intersection of all maximal subalgebras of $A$  if those 
exist and as $A$ otherwise.
G.L.~Carns and C.~-Y.~Chao \cite{2} showed that if $A$ is a nilpotent algebra over a
field $K$,  then $\Phi(A)=A^2$. It follows that 
\itemitem{}(III) If $B$ is a filtered multiplicative $K$-basis  of $A$
and if $B\setminus\{1\}\subseteq rad(A)$, then all elements of 
$B\setminus rad(A)^2$ are generators of $A$ over  $K$.
 
Now  let $G$ be a group. For $a,b\in G$ we define
${}^ba=bab\m1$ and  $[a,b]=aba\m1b\m1$. 
The ideal 
$$
I_K(G)=\big\{\,\,\, \sum_{g\in G}\alpha_gg\in KG \,\,\,\, \mid\,\,\,\,  
\sum_{g\in G}{\alpha_g}=0\,\,\, \big\}
$$ 
is called the augmentation ideal of $KG$. Then the following subgroup
$$
{\frak D}_n(G)=\big\{\,\,\, g\in G \,\,\,\, \mid \,\,\,\,  
g-1\in I^n_K(G)\,\,\, \big\}.
$$
is called the $n$-th dimensional subgroup of $KG$.

\subhead
3. Results
\endsubhead
By theorem 3.11.2 of \cite{3} every metacyclic $p$-group has the following
presentation 
$$
G=\gp{a,b\,\,\,\, \mid \,\,\,\,  a^{p^n}=1, \,\,\,\,   b^{p^m}=a^{p^t}, \,\,\,\,   {}^ba=a^r},
$$
where $t\geq0$, $r^{p^m}\equiv 1\pmod{p^n}$ and $p^t(r-1)\equiv 0\pmod
{p^{n}}$.
Therefore, every element of $G$ can be written as $g=a^{i}b^j$, where 
$0\leq i< p^n$ and  $0\leq j<p^{m-t}$. Using the identity 
$$
(xy-1)=(x-1)(y-1)+(x-1)+(y-1), \tag1
$$ 
we obtain that every element of the   augmentation ideal $I_K(G)$ is a sum of
elements of  the form $(a-1)^{k}(b-1)^l$, where  $0\leq k< p^n$,  $0\leq
l<p^{m-t}$ and $k+l\geq 1$.

\proclaim
{Theorem}
Let $G$ be a finite metacyclic $p$-group and  $K$ be a field of characteristic
$p$.  
Then the group algebra $KG$ possesses a filtered multiplicative $K$-basis  if
and only if $p=2$ and exactly one of the
following conditions holds:
\itemitem{1.} $G$ is a dihedral group;
\itemitem{2.} $K$  contains  a   primitive
cube root of the unity and $G$ is a quaternion group of order $8$.
\endproclaim

\demo{Proof} 
Clearly, $I_K(G)$ is a radical of $KG$. Suppose that $\{1,B\}$ is a filtered
multiplicative $K$-basis  of $KG$. Then $B$ is a filtered multiplicative
$K$-basis  of $I_K(G)$.  Obviously, $(a-1)^{i}(b-1)^j\in I^2_K(G)$ if $i+j\geq
2$ and $a-1, b-1$ are generators  of $I_K(G)$ over $K$. By Jennings
theory \cite{4}, $(a-1)+I^2_K(G)$ and $(b-1)+I^2_K(G)$ form a  $K$-basis of
$I_K(G)/I^2_K(G)$. Therefore, by property (III), the subset $B\setminus B^2$
consists of two elements, which we  denote  $u$ and $v$. Thus $K[u,v]=I_K(G)$
and  
$$
\cases
u\equiv \alpha_1(a-1)+\alpha_2(b-1)\pmod{I^2_K(G)};\\
v\equiv \beta_1 (a-1)+\beta_2 (b-1)\pmod{I^2_K(G)},
\endcases\tag2
$$
where $\alpha_i,\beta_i\in K$ and
$\Delta=\alpha_1\beta_2-\alpha_2\beta_1\not=0$.

Clearly, $c=[b,a]\in {\frak D}_2(G)$ and $c-1\in I^2_K(G)$. By a simple
calculation we get  
$$
uv\equiv\alpha_1\beta_1(a-1)^2+\alpha_2\beta_2(b-1)^2+
$$$$
+(\alpha_1\beta_2+\alpha_2\beta_1)(a-1)(b-1)+ \alpha_2\beta_1(c-1) \pmod{I^3_K(G)},\tag3
$$$$
vu\equiv\alpha_1\beta_1(a-1)^2+\alpha_2\beta_2(b-1)^2+
$$$$
+(\alpha_1\beta_2+\alpha_2\beta_1)(a-1)(b-1)+ \alpha_1\beta_2(c-1) \pmod{I^3_K(G)},\tag4
$$$$
u^2\equiv\alpha_1^2(a-1)^2+\alpha_2^2(b-1)^2+
2\alpha_1\alpha_2(a-1)(b-1)+ \alpha_1\alpha_2(c-1) \pmod{I^3_K(G)},\tag5
$$$$
v^2\equiv\beta_1^2(a-1)^2+\beta_2^2(b-1)^2+
2\beta_1\beta_2(a-1)(b-1)+ \beta_1\beta_2(c-1) \pmod{I^3_K(G)}.\tag6
$$

We consider the case  when $c-1\in I^3_K(G)$. Then by (3) and (4) we have  $uv\equiv
vu \pmod {I^3_G(K)}$.  Moreover,  $uv, vu\not \in I^3_K(G)$. Indeed, 
if $uv$ or  $vu \in I^3_K(G)$  then by (3) or (4) we obtain   
$\alpha_1\beta_1= \alpha_2\beta_2=\alpha_1\beta_2+\alpha_2\beta_1=0$ and
$\Delta=0$, which is impossible. 
Therefore, $uv,vu\not\in I^3_K(G)$ and  $uv\equiv vu \pmod {I^3_K(G)}$ and by
property (II) of the filtered multiplicative $K$-basis  we conclude that $uv=vu$
and $I_K(G)$ is a commutative algebra, which is contradiction. 

In the rest of the proof we assume that $c-1\not\in I^3_K(G)$. It is well-known
that  for all nonabelian $p$-groups  the factor group $G/G'$ is not cyclic (see \cite{3},
Theorem 3.7.1). Thus 
$r-1$ is divisible by $p$ and $r-1=ps$ for some $s$. Then $c-1=(a^s-1)^p\in
I^3_K(G)$ for $p>2$ and also for $p=2$  if $s$ is even. We have established that
$s$ is odd and $G$ is a $2$-group with the following defining presentation:
either 
$$
G=\gp{\,\,\, a,b\,\,\,\, \mid  \,\,\,\,  a^{2^n}=1, \,\,\,\,   {}{}b^{2^m}=1, \,\,\,\,   
{}^ba=a^r\,\,\, }, \tag7
$$
where $r^{2^m}\equiv 1 \pmod{2^n}$ and $r\not\equiv 1\pmod 4$,
or 
$$
G=\gp{\,\,\,  a,b\,\,\,\, \mid  \,\,\,\,  a^{2^n}=1, \,\,\,\,   
b^{2^m}=a^{2^{n-1}}, \,\,\,\,   {}^ba=a^r\,\,\, },\tag8
$$
where $r^{2^m}\equiv 1 \pmod{2^n}$, $2^{n-1}(r-1)\equiv 0 \pmod {2^n}$ and $4$
does  not divide $r-1$.

Suppose that $G$ has the defining presentation (7) and $b^2=1$. Since $r-1=2s$
and $(s,2)=1$, from $r^2\equiv 1\pmod {2^n}$  it follows that   $s=-1$ or
$s=-1+2^{n-2}$  for  $n\geq 3$. Then by (1) we have   
$
c+1=a^{2s}+1\equiv (1+a)^2 \pmod {I^3_K(G)}
$  
and it follows from (3)--(6) that 
$$
\cases
uv\equiv \beta_1(\alpha_1+\alpha_2)(1+a)^2+\Delta(1+a)(1+b)&\pmod{I^3_K(G)};\cr
vu\equiv \alpha_1(\beta_1+\beta_2)(1+a)^2+\Delta(1+a)(1+b)&\pmod{I^3_K(G)};\cr
u^2\equiv \alpha_1(\alpha_1+\alpha_2)(1+a)^2&\pmod{I^3_K(G)};\cr
v^2\equiv \beta_1(\beta_1+\beta_2)(1+a)^2&\pmod{I^3_K(G)}.\cr
\endcases \tag9
$$

Clearly  $uv,vu\not\in I^3_K(G)$ and by $\Delta\not=0$ we have that
$uv\not\equiv vu \pmod{I^3_K(G)}$. Since the $K$-dimension of
$I^2_K(G)/I^3_K(G)$ equals $2$, the elements $uv+I^3_K(G)$ and  $vu+ I^3_K(G)$
form a $K$-basis of $I^2_K(G)/I^3_K(G)$ and $u^2,v^2\in I^3_K(G)$. We conclude
that  
$
\alpha_1(\alpha_1+\alpha_2)=0$  and  $\beta_1(\beta_1+\beta_2)=0
$,  
whence it  follows that $u= \alpha(a+b)$ and
$v=\beta(1+b)$. Clearly we can set  $\alpha=\beta=1$.

Let $G=\gp{ \,\,\, a,b\,\,\,\, \mid  \,\,\,\,  
a^{2^n}=1,  \,\,\,\,  {}{}b^{2}=1,  \,\,\,\,  {}^ba=a^{-1}
\,\,\, }$ with $n\geq 2$   be
a dihedral group of order $2^{n+1}$. 
We shall prove   by induction
in $i$ that  $u^{i}$ can be written as 
$$
(1+a)^{2i-1}+(1+a)^{2i-2}(1+b)+\beta_1(1+a)^{2i}+ 
$$$$
+\beta_2(1+a)^{2i-1}(1+b) \pmod{I^{2i+1}_K(G)},\tag10
$$
where $\beta_1=\beta_2=1$ if $i$ is even and 
$\beta_1=\beta_2=0$ otherwise.

Base of induction:
It is easy to see that this is true for $i=1,2$, and the induction  
step follows by,  
$$(1+b)(1+a)
\equiv (1+a)(1+b)+(1+a)^2(1+b)
+(1+a)^2+
$$$$
+(1+a)^3(1+b)+(1+a)^3+(1+a)^4 \pmod{I^5_K(G)}
$$  
and 
$$
u^{i}u
\equiv(\beta_1+\beta_2+1)[(1+a)^{2i+1}+(1+a)^{2i}(1+b)]+
$$$$
+(1+\beta_2)[(1+a)^{2i+2}+(1+a)^{2i+1}(1+b)]  \equiv u^{i+1}
\pmod{I^{2i+3}_K(G)}. 
$$
Hence   (10) holds.

Using (10),  we obtain that
$$
u^{i}\equiv (1+a)^{{2i-1}}+(1+a)^{{2i-2}}(1+b)\pmod{I^{2i}_K(G)},
$$$$ 
vu^{i}\equiv  (1+a)^{{2i}}+(1+a)^{{2i-1}}(1+b)\pmod{I^{2i+1}_K(G)},
$$$$
u^{j}v\equiv (1+a)^{{2j-1}}(1+b)\pmod{I^{2j+1}_K(G)},
$$$$
vu^{j}v\equiv (1+a)^{{2j}}(1+b)\pmod{I^{2j+2}_K(G)}, 
$$
where $i=1,\ldots, 2^{n-1}$ and $j=1,\ldots, 2^{n-1}-1$. 

Clearly, the factor algebra $I^{t}_K(G)/I^{t+1}_K(G)$ has the following basis:
$(a+1)^t+I^{t+1}_K(G)$ and $(a+1)^{t-1}(b-1)+I^{t+1}_K(G)$.

First, let $t=2k+1$, where $k=1,\ldots, 2^{n-2}-1$. Then we have 
$$
u^{k+1}\equiv (1+a)^{2k+1}+(1+a)^{2k}(1+b) \pmod{I^{t+1}_K(G)},
$$$$
vu^{k}v\equiv              (1+a)^{2k}(1+b) \pmod{I^{t+1}_K(G)}
$$
and it follows that $u^{k+1}$ and $vu^{k}v$ are linearly independent by modulo
$I^{t+1}_K(G)$.

Now, let $t=2k$, where $k=1,\ldots,  2^{n-2}-1$. Then we have 
$$
vu^{k}\equiv (1+a)^{2k}+(1+a)^{2k-1}(1+b) \pmod{I^{t+1}_K(G)},
$$$$
u^{k}v\equiv              (1+a)^{2k-1}(1+b) \pmod{I^{t+1}_K(G)}
$$
and, as before,  $vu^{k}$ and $u^{k}v$ are linearly independent by modulo
$I^{t+1}_K(G)$.

Therefore  the matrix of decomposition is unitriangle and 
$$
\{ 1,v, u^{i}, vu^{i}, u^{j}v, vu^{j}v \,\,\,\, \mid  
i=1,\ldots, 2^{n-1} \text{ and } j=1,\ldots, 2^{n-1}-1 \}
$$  
form a filtered multiplicative $K$-basis of $KG$.

Now let $G=\gp{a,b\,\,\,\, \mid  \,\,\,\,  a^{2^n}=b^2=1,  \,\,\,\,  {}^ba=a^{-1+2^{n-1}}}$ with $n\geq 3$  be a
semidihedral group and set $u=a+b$, $v=1+b$. An easy calculation gives 
$1+a\m1=\sum_{i=1}^{2^n-1}(1+a)^i$ and 
$$
u^2=\sum_{i=2}^{2^{n-1}-1}(1+a)^{i}(1+b)+ \sum_{j=3}^{2^{n-1}-1}(1+a)^j, 
$$ 
$$
uvu= \sum_{i=2}^{2^{n-1}}(1+a)^{i}(1+b)+\sum_{j=3}^{2^{n-1}}(1+a)^j.
$$ 
Therefore $u^2\equiv uvu\pmod{I_K^{2^{n-1}}(G)}$, but $u^2,uvu\not\in
I_K^4(G)$ and 
$$
u^2- uvu=(1+a)^{2^{n-1}}(1+b)+(1+a)^{2^{n-1}}\not=0,  
$$ 
which  contradicts  property (II). 

Suppose that $G$ has the defining presentation (7) with $m>1$ or (8) with $m>1$. 
By (1) we have 
$$
(1+b)(1+a)
\equiv(1+a)(1+b)+(1+a^s)^2 \equiv
(1+a)(1+b)+(1+a)^2\pmod{I^3_K(G)}
$$ 
and  it follows from  (3)--(6) that 
$$
\cases
uv\equiv \beta_1(\alpha_1+\alpha_2)(1+a)^2+\Delta(1+a)(1+b)+
\alpha_2\beta_2(1+b)^2&\pmod{I^3_K(G)};\cr 
vu\equiv \alpha_1(\beta_1+\beta_2)(1+a)^2+\Delta(1+a)(1+b)+
\alpha_2\beta_2(1+b)^2&\pmod{I^3_K(G)};\cr 
u^2\equiv \alpha_1(\alpha_1+\alpha_2)(1+a)^2+
\alpha_2^2(1+b)^2&\pmod{I^3_K(G)};\cr 
v^2\equiv \beta_1(\beta_1+\beta_2)(1+a)^2+
\beta_2^2(1+b)^2&\pmod{I^3_K(G)}.\cr 
\endcases  \tag11
$$

It is easy to see  that  $uv,vu\not\in I^3_K(G)$. Using  the fact that 
$\Delta\not=0$,  we establish  that
$uv\not\equiv vu \pmod{I^3_K(G)}$. Therefore  $uv+I^3_K(G)$ and  $vu+
I^3_K(G)$ 
are  $K$-linearly independent. It is easily verified that 
$u^2+I^3_K(G)$ and  $v^2+ I^3_K(G)$ are nonzero elements of  
 $I^2_K(G)/I^3_K(G)$ and $uv\not\equiv v^2$, $vu\not\equiv v^2$, 
$uv\not\equiv u^2$, $vu\not\equiv u^2$.  Since the $K$-dimension of
$I^2_K(G)/I^3_K(G)$  equals  $3$, we have $u^2\equiv v^2\pmod{I^3_K(G)}$ 
 and by property (II) of  the filtered multiplicative $K$-basis,
 $u^2=v^2$. From $u^2\equiv v^2 \pmod{I^3_K(G)}$ we  
obtain $\alpha_2^2=\beta_2^2$ and
$\alpha_1(\alpha_1+\alpha_2)=\beta_1(\beta_1+\beta_2)$.  
By $\Delta\not=0$ we have  $\alpha_2=\beta_2\not=0$,  
whence  the equation   
$\alpha_1^2+\beta_2\alpha_1+\beta_1(\beta_1+\beta_2)=0$  
has a solution $\alpha_1=\beta_1+\beta_2$  whence   $\Delta=\beta_2^2\not=0$.  
Thus we observe that
$u=(1+\lambda)a+b+\lambda$ and $v=\lambda{a}+b+\lambda+1$, where
$\lambda=\frac{\beta_1}{\beta_2}$.  
Then,   keeping the equality $u^2=v^2$,  we conclude that  
$1+a^2+ab+ba=0$, which is impossible.

Suppose that $G$ has the defining presentation (8) with $m=1$.  As we obtained
before,  either $r=-1$ or $r=-1+2^{n-1}$. By (1) we have 
$$
(1+b)(1+a)\equiv (1+a)(1+b)+(1+a)^2 \pmod{I^3_K(G)}
$$
and we can write the elements $u,v$ in the  form (2). It follows from (3)--(6)
that  (11) hold by modulo $I^3_K(G)$.

We shall  consider two  cases depending on the values of $r$ and $m$ in (8).

Case 1. Let $G$ be a quaternion group of order $8$. Then by (11) we have 
$$
\cases
uv\equiv
(\alpha_1\beta_1+\alpha_2\beta_1+\alpha_2\beta_2)(1+a)^2+\Delta(1+a)(1+b) \pmod{I^3_K(G)};\\ 
vu\equiv 
(\alpha_1\beta_1+\alpha_1\beta_2+\alpha_2\beta_2)(1+a)^2+\Delta(1+a)(1+b) \pmod{I^3_K(G)};\\ 
u^2\equiv 
(\alpha_1^2+\alpha_1\alpha_2+\alpha_2^2)(1+a)^2 \pmod{I^3_K(G)};\\  
v^2\equiv 
(\beta_1^2+\beta_1\beta_2+\beta_2^2)(1+a)^2 \pmod{I^3_K(G)}.\\  
\endcases
$$

Since the $K$-dimension of $I^j_K(G)/I^{j+1}_K(G)$ ($j=1,\ldots, 4$) equals
$2$ and $uv\not\equiv vu\pmod{I^3_K(G)}$, we have
$\alpha_1^2+\alpha_1\alpha_2+\alpha_2^2$ and
$\beta_1^2+\beta_1\beta_2+\beta_2^2=0$. Using the fact that $\Delta\not=0$,
we establish  $\frac{\alpha_1}{\alpha_2}=\omega$,
$\frac{\beta_2}{\beta_1}=\omega^2$. Thus we observethat
$u=\omega(1+a)+(1+b)$ and  
$v=(1+a)+\omega^2(1+b)$, where $\omega$ is a primitive cube root of the
unity. 

A simple calculation by modulos  $I^4_K(G)$, $I^5_K(G)$ shows that 
$$
\{1,u,v,  uv, vu, uvu, vuv, uvuv \}
$$ 
is  a filtered multiplicative $K$-basis for $KG$.

Case 2.  Let $G$ has a presentation 
$$
\gp{a,b\,\,\,\, \mid  \,\,\,\,  a^{2^n}=1,  \,\,\,\,  b^2=a^{2^{n-1}},
{}^ba=a^r}\tag12 
$$ 
with $n>2$. Then by (11) we have 
$$
\cases
uv\equiv
(\alpha_1+\alpha_2)\beta_1(1+a)^2+\Delta(1+a)(1+b) \pmod{I^3_K(G)};\\ 
vu\equiv 
\alpha_1(\beta_1+\beta_2)(1+a)^2+\Delta(1+a)(1+b) \pmod{I^3_K(G)};\\ 
u^2\equiv 
\alpha_1(\alpha_1+\alpha_2)(1+a)^2 \pmod{I^3_K(G)};\\  
v^2\equiv 
\beta_1(\beta_1+\beta_2)(1+a)^2 \pmod{I^3_K(G)}.\\  
\endcases
$$

Since the $K$-dimension of $I^2_K(G)/I^3_K(G)$  equals $2$ and $\Delta\not=0$,
we have either  $\alpha_1=\alpha_2\not=0$ and $\beta_1=0$ or  $\alpha_1=0$ and
$\beta_1=\beta_2\not=0$.  The second case is  similar to  first. Therefore,
we can put $u=(1+a)+(1+b)$, $v=1+b$.

Case 2.1. Let $r=-1$ in (12). Then $G$ is a generalized quaternion group. An easy
calculation gives 
$$
(1+b)(1+a)=\sum_{j=1}^{2^{n}-1}(1+a)^j(1+b)+\sum_{j=1}^{2^n-1}(1+a)^{j+1}
$$
and 
$$
u^2=\sum_{j=1}^{2^{n}-1}(1+a)^{j+1}(1+b)+\sum_{j=1}^{2^{n}-1}(1+a)^{j+2}+(1+a)^{2^{n-1}},
$$ 
$$
uvu=\sum_{j=1}^{2^{n}-1}(1+a)^{j+1}(1+b)+\sum_{j=1}^{2^{n}-1}(1+a)^{j+2}+(1+a)^{2^{n-1}}(1+b).
$$

Therefore, $u^2\equiv uvu\pmod{I^4_K(G)}$, but $u^2,uvu\not\in I^4_K(G)$ and 
$$
u^2-uvu=(1+a)^{2^{n-1}}
+(1+a)^{2^{n-1}}(1+b)\not=0,
$$
which  contradicts property (II).

Case 2.2. Let $G=\gp{a,b\,\,\,\, \mid  \,\,\,\,  a^{2^n}=1,  \,\,\,\,  b^2=a^{2^{n-1}}, \,\,\,\,  
{}^ba=a^{-1+2^{n-1}}}$. It is easy to see that $(ab)^2=a^{2^{n-1}}b^2=1$ and 
$$
G\cong \gp{\,\,\,  a,ab\,\,\,\, \mid  \,\,\,\,  a^{2^n}=1,  \,\,\,\,  (ab)^2=1,  \,\,\,\, 
 {}^{ab}a=a^{-1+2^{n-1}}},
$$ 
which is a semidihedral group and, as we  saw  before,  $KG$ has no  filtered
multiplicative $K$-basis.

Thus our theorem is proved. 
\enddemo
\subhead
4. Example
\endsubhead
Now we give an  example of a nonmetacyclic $2$-group with  a filtered multiplicative
basis. 

Let $G=\gp{\,\,\, a,b\,\,\,\, \mid  \,\,\,\,  a^4=b^4=1,  \,\,\,\,  {}^ba=b^2a^3,  \,\,\,\,  {}^{a}b=a^2b^3, \,\,\,\,  
[a^2,b]=[b^2,a]=1\,\,\, }$, a group of order $16$,  and let  $K$ be a field of characteristic $2$.  
Then elements 
$$
\{1, u, v, uv, vu, v^2, uvu,uv^2, vuv,v^3, uvuv, uv^3, vuv^2,
uvuv^2, vuv^3, uvuv^3\,\,\,\, \mid 
$$$$
u=a+b, v=\mu_1a+\mu_2b+(\mu_1+\mu_2)\text{ and } \mu_1,\mu_2\in K,\text{ and } \mu_1\not=\mu_2\}
$$
form a filtered multiplicative $K$-basis for $KG$.

Indeed, by  (1) we have 
$$
(1+b)(1+a)\equiv (1+a)(1+b)+(1+a)^2+(1+b)^2 \pmod{I^3_K(G)}
$$
and $u,v$ be  can writen in the  form (2).    

By a simple calculation modulo $I^3_K(G)$ we have
$$
\cases
uv\equiv (\alpha_1+\alpha_2)\beta_1(1+a)^2+\Delta(1+a)(1+b)+\alpha_2(\beta_1+\beta_2)(1+b)^2;\\
vu\equiv \alpha_1(\beta_1+\beta_2)(1+a)^2+\Delta(1+a)(1+b)+\beta_2(\alpha_1+\alpha_2)(1+b)^2;\\
u^2\equiv (\alpha_1+\alpha_2)\alpha_1(1+a)^2+\alpha_2(\alpha_1+\alpha_2)(1+b)^2;\\
v^2\equiv (\beta_1+\beta_2)[\beta_1(1+a)^2+\beta_2(1+b)^2].\\
\endcases\tag13
$$
It is easy to see that $K$-dimension of $I^2_K(G)/I^3_K(G)$ equals $3$ and
\newline
$uv\not\equiv vu\pmod{I^3_K(G)}$,  
$uv\not\equiv u^2\pmod{I^3_K(G)}$,
$uv\not\equiv v^2\pmod{I^3_K(G)}$,
\newline
$vu\not\equiv u^2\pmod{I^3_K(G)}$,
$vu\not\equiv v^2\pmod{I^3_K(G)}$.

 We have the following two cases.

First let $u^2\equiv v^2\not\equiv 0\pmod{I^3_K(G)}$. Then by (13) we have 
$\alpha_1^2+ \alpha_1\alpha_2=\beta_1^2+\beta_1\beta_2$ 
and 
$\alpha_2^2+\alpha_1\alpha_2=\beta_2^2+\beta_1\beta_2$. It follows that 
$(\alpha_1+\alpha_2)^2=(\beta_1+\beta_2)^2$ 
and $\alpha_1+\alpha_2=\beta_1+\beta_2$. Then by $u^2\equiv v^2\not\equiv
0\pmod{I^3_K(G)}$ we have  $\Delta=0$,  which is impossible.

Now let $u^2\equiv 0\pmod{I^3_K(G)}$ or $v^2\equiv 0\pmod{I^3_K(G)}$.  It is
easy to see that the second case is symmetric to 
the first one, so we consider only the first case.  
Then $\alpha_1=\alpha_2\not=0$ and by (13) we have 

$$
\cases
uv \equiv\lambda[(1+a)(1+b)+(1+b)^2]\pmod{I^3_K(G)};\\
vu \equiv\lambda[(1+a)^2+(1+a)(1+b)]\pmod{I^3_K(G)};\\
v^2\equiv\lambda[\beta_1(1+a)^2+\beta_2(1+b)^2]\pmod{I^3_K(G)},\\
\endcases
$$
where $\lambda=\beta_1+\beta_2\not=0$. 
By a simple calculation  modulo $I^4_K(G)$ we obtain  

$$
\cases
uvu \equiv \lambda[   (1+a)^3 +      (1+a)^2(1+b) +    (1+a)(1+b)^2  +    (1+b)^3];\\
uv^2 \equiv\lambda[\beta_1(1+a)^3 + \beta_1(1+a)^2(1+b) +\beta_2(1+a)(1+b)^2  +\beta_2(1+b)^3];\\
vuv \equiv \lambda^2[                  (1+a)^2(1+b) +    (1+a)(1+b)^2       ];\\
v^3\equiv  \lambda[\beta_1^2(1+a)^3 + \beta_1\beta_2(1+a)^2(1+b) +\beta_1\beta_2(1+a)(1+b)^2  +\beta_2^2(1+b)^3]\\
\endcases
$$
 and  modulo $I^5_K(G)$ 
$$
\cases
uvuv  \equiv \lambda^2[(1+a)^3(1+b)+(1+a)^2(1+b)^2+(1+a)(1+b)^3] ;\\
uv^3  \equiv \lambda^2[\beta_1(1+a)^3(1+b)+\beta_1(1+a)^2(1+b)^2+\beta_2(1+a)(1+b)^3];\\
vuv^2 \equiv \lambda^2[\beta_1(1+a)^3(1+b)+\beta_2(1+a)^2(1+b)^2+\beta_2(1+a)(1+b)^3].\\
\endcases
$$
Similarly
$$
\cases
uvuv^2    \equiv \lambda^2[(1+a)^3(1+b)^2+(1+a)^2(1+b)^3]\pmod{I^6_K(G};\\
vuv^3    \equiv \lambda^2[\beta_1(1+a)^3(1+b)^2+\beta_2(1+a)^2(1+b)^3]\pmod{I^6_K(G};
\endcases
$$
and 
$uvuv^3\equiv \lambda^3(1+a)^3(1+b)^3\pmod{I^7_K(G)}$.

Since the number of elements  modulo  $I^j_K(G)$, ($j=2,\ldots, 6$) equals
the numbers of the $K$-dimension of $I^j_K(G)/I^{j+1}_K(G)$, we conclude  that the
elements $\{1, u, v, uv, vu, v^2, uvu,uv^2, vuv,v^3, uvuv, uv^3, vuv^2,
uvuv^2, vuv^3, uvuv^3\}$ form a filtered multiplicative $K$-basis for $KG$.

\enddocument